\documentclass[12pt,a4paper]{amsart}

\usepackage{amsfonts, amsmath, amssymb, amsthm, amscd, hyperref}

\usepackage{a4wide}

\newtheorem{thm}{Theorem}

\newtheorem{con}{Conjecture}

\theoremstyle{definition}
\newtheorem{defn}{Definition}

\theoremstyle{remark}

\renewcommand{\int}{\mathop{\rm int}}

\renewcommand{\epsilon}{\varepsilon}

\begin{document}
\title{Inscribing a regular crosspolytope}
\author{R.N.~Karasev}

\thanks{The research of R.N.~Karasev was supported by the President's of Russian Federation grant MK-1005.2008.1, and partially supported by the Dynasty Foundation.}

\email{r\_n\_karasev@mail.ru}
\address{
Roman Karasev, Dept. of Mathematics, Moscow Institute of Physics
and Technology, Institutskiy per. 9, Dolgoprudny, Russia 141700}

\begin{abstract}
In this paper we prove that it is possible to inscribe a regular crosspolytope (multidimensional octahedron) into a smooth convex body in $\mathbb R^d$, where $d$ is an odd prime power. Some generalizations of this statement are also proved.
\end{abstract}


\maketitle

\section{Introduction}

The problems of inscribing or outscribing a polytope of a given family to some smooth or convex line or surface have a long history. In~\cite{shn1944} it was proved that a square can be inscribed into any simple smooth closed curve in the plane (Schnirelmann's theorem). In~\cite{kak1942} it was proved that a cube can be outscribed about any convex body in $\mathbb R^3$. In~\cite{gro1969} it was shown that if a smooth closed hypersurface $H\subset \mathbb R^d$ bounds a region with nonzero Euler characteristic, then $H$ contains vertices of a simplex, positively homothetic to any prescribed simplex. 

The questions of inscribing and outsribing are widely discussed in the books~\cite{gru1971,kleewa1996}. More recent results and references can be found in~\cite{mak2003,mak2006}. 

The inscribing and outscribing problems are closely connected to the Knaster problem~\cite{kna1947} on the level surfaces of function on $S^{d-1}$.

\begin{con}[Knaster's problem for functions]
\label{knaster}
Let $S^{d-1}$ be a unit sphere in $\mathbb R^d$. Suppose that the points $x_1, \ldots, x_d\in S^{d-1}$ and a continuous function $f: S^{d-1}\to \mathbb R$ are given. Then there exists a rotation with positive determinant $\rho\in SO(d)$ such that
$$
f(\rho(x_1)) = f(\rho(x_2)) = \dots = f(\rho(x_d)).
$$
\end{con}

In~\cite{kasha2003,hiri2005} some counterexamples to the Knaster problem were found, with some special functions and point sets, the dimension $d$ was required to be large. Still for some particular point sets $\{x_1, \ldots, x_d\}$ the conjecture turns out to be true. In Section~\ref{obstruction} we use a positive solution to the Knaster problem, that is formulated in terms of the Euler class of some vector bundle.

In the book~\cite{kleewa1996} the problem of inscribing a regular octahedron into a convex body in $\mathbb R^3$ is formulated (Problem 11.5). For the case of $d=3$ and a smooth convex body, this problem was solved in~\cite{mak2003}. Here we prove a similar statement for larger dimensions.

\begin{defn}
Let an \emph{orthogonal cross} in $\mathbb R^d$ be $d$ mutually orthogonal straight lines, passing through a single point $o$, the latter is called \emph{the center of the cross}.
\end{defn}

\begin{defn}
Let $(e_1,\ldots, e_d)$ be some linear base in $\mathbb R^d$. 
The convex hull of the points
$$
e_1, -e_1, e_2, -e_2,\ldots, e_d, -e_d\in\mathbb R^d
$$
and all its similar images (under transforms with positive determinant) are called a \emph{crosspolytope} in $\mathbb R^d$. If the base is orthonormal, we call the crosspolytope \emph{regular}.
\end{defn}

\begin{defn}
A convex body $K\subset\mathbb R^d$ is called \emph{non-angular}, if none of the points $o\in\partial K$ is a center of some orthogonal cross $C$ such that $C\cap\int K=\emptyset$.
\end{defn}

It is clear that smooth convex bodies are non-angular. Now we can formulate the main results.

\begin{thm}
\label{inscrcross}
Let $K\subset\mathbb R^d$ be a non-angular convex body, let $d$ be an odd prime power. Then there exists a regular crosspolytope $C\subset\mathbb R^d$ such that all its vertices are on $\partial K$.
\end{thm}

The same proof is used to establish the following fact.

\begin{thm}
\label{insrccross-p}
Let $H\subset\mathbb R^d$ be an image of some smooth embedding of $S^{d-1}$, where $d=p^k$ is an odd prime power. Let $(e_1,\ldots, e_d)$ be some linear base in $\mathbb R^d$, let $C$ be the convex hull of $(\pm e_1,\ldots, \pm e_d)$. 

Suppose that the group $G=(Z_p)^k$ acts transitively on the vectors $(e_1,\ldots, e_d)$ and this action is a restriction of some (special) orthogonal action of $G$ on $\mathbb R^d$. Then there exists a crosspolytope $C'\subset\mathbb R^d$, similar to $C$, with all its vertices lying on $H$.
\end{thm}

It should be mentioned that the possibility to inscribe a regular crosspolytope into a centrally symmetric convex body follows from the positive solution of the Knaster problem for the orthonormal base in~\cite{yayu1950}, which is true for every dimension $d$.

Theorem~\ref{insrccross-p} for a regular crosspolytope in arbitrary dimension was announced in~\cite{gug1965}, but the sketch of the proof was based on some two-dimensional reasoning. Though, the two-dimensional case (Schnirelmann's theorem) was proved in~\cite{gug1965} incorrectly, the main lemma on the continuous dependence of the inscribed square on the smooth curve is false and it has simple counterexamples. 
In the plane~\cite{shn1944} the parity of the number of inscribed squares remains constant under the deformations of the curve. In the space of higher dimension there is an infinite variety of inscribed crosspolytopes, and the counting and parity argument does not work.

In the case $d=2$ the proof in this paper cannot be applied. In this case some larger group (at least $Z_4$) should be considered as the symmetry group. 

\section{Reduction to a topological fact}

Denote $\{1,2,\ldots,n\} = [n]$.

First we outline the proof of Theorem~\ref{inscrcross} for strictly convex non-angular $K$, which also proves Theorem~\ref{insrccross-p} for the case, when $H$ is a boundary of some smooth strictly convex body $K$.

Let $R$ be the space of positively oriented orthonormal frames in $\mathbb R^d$. In the case of Theorem~\ref{insrccross-p} let $R$ be the space of frames, that can be obtained from the given frame $(e_1,\ldots, e_d)$ by a special orthogonal transform. In both cases $R$ is naturally identified with $SO(d)$. Note that in Theorem~\ref{insrccross-p} the lengths of the vectors are equal, so we consider them to be unit vectors.

Consider a continuous map $g$ from $K\times R$ to the linear space $V=\mathbb R^{2d}$, defined as follows. For a point $p\in K$, a frame $(e_1,\ldots,e_d)\in R$ and any $i\in [d]$ put
$$
a_i = \max \{a : p + ae_i\in K\},\quad b_i = \max \{a : p - ae_i\in K\}.
$$
It is quite clear that for strictly convex $K$ this map is continuous.

Consider the coordinates $s_i = a_i + b_i, t_i = a_i - b_i$ in $V$, and consider the one-dimensional subspace $L\subset V$, given by
$$
t_1=\dots=t_d = 0,\quad s_1=\dots=s_d.
$$
In the quotient space $V/L$, let the $d$-dimensional linear hull of $\{t_1,\ldots, t_d\}$ be $U$, and the $d-1$-dimensional linear hull of $\{s_1,\ldots, s_d\}$ be $W$. Now we have to prove that the quotient map $f: K\times R\to U\oplus W$ maps some pair $(p, r)\in \int K\times R$ to zero. 

\section{Calculating the obstruction}
\label{obstruction}

In the sequel we consider cohomology $\mod p$, the coefficients are omitted in the notation. The necessary facts on algebraic topology and vector bundles can be found in~\cite{hsiang1975,milsta1974,mishch1998}.

Note that if $K$ is non-angular (or smooth), the map $f$ is nonzero at any pair $(p, r)\in \partial K\times R$. In case of Theorem~\ref{inscrcross} we can consider some free transitive action of $G=(Z_p)^k$ (where $d=p^k$) on $R$ by permuting the vectors. In case of Theorem~\ref{insrccross-p} the action of $G$ on $R$ is already given. Correspondingly, $G$ acts on $U$ and $W$ by permuting the respective coordinates $t_i, s_i$.

Thus the map $f$ is an equivariant section of a $G$-bundle 
$$
K\times R\times U\times W\to K\times R
$$ 
over the $G$-space $K\times R$ and the obstruction for $f$ being nonzero over the entire $K\times R$ is the relative Euler class (see~\cite{ker1957,kar2008dcpt}), which resides in the equivariant cohomology $H_G^{2d-1}(K\times R, \partial K\times R)$.

The section $f$ splits as a direct sum of $s_U\oplus s_W$ for the corresponding $G$-bundles. Under this splitting the section $s_U$ is nonzero over $\partial K\times R$ and gives some relative Euler class in $e(s_U)\in H_G^d(K\times R, \partial K\times R)$. The section $s_W$ has some Euler class $e(s_W)\in H_G^{d-1} (K\times R) = H_G^{d-1}(R)$. The Euler class of $f$ is therefore $e(f) = e(s_U)e(s_W)$ by the multiplicative rule for the Euler class.

It is clear from the K\"unneth formula that the cohomology $H_G^*(K\times R, \partial K\times R)$ equals the tensor product $H^*(K,\partial K)\otimes H_G^*(R)$, in this particular case it is $u \times H_G^*(R)$, where $u$ is the $d$-dimensional generator of $H^*(K,\partial K)$. 

Let us find $e(s_U)$ first. Note that if we deform $K$, leaving it non-angular, this class remains the same. Thus it is sufficient to consider $K$ equal to the unit ball $B$. In this case we fix the frame $r$ in the pair $(p,r)\in B\times R$ and note that the section $s_U$ has the only non-singular zero in the center of $B\times\{r\}$. Hence the image of $e(s_U)$ under the natural map  $H_G^d(B\times R, \partial B\times R)\to H^d(B,\partial B)$ gives the generator of $H^d(B,\partial B)$, hence $e(s_U) = u\times 1 \in H_G^*(K\times R, \partial K\times R)$.

Now note that the Euler class $e(s_W)$ coincides with the obstruction, which is used in establishing a particular case of the Knaster problem in~\cite{vol1992}. By definition $W=\mathbb R[G]/\mathbb R$ as a representation of $G$, this representation has no trivial irreducible summands. Therefore (see~\cite[Ch.~III~\S~1]{hsiang1975}) the Euler class of this representation in the cohomology of the classifying space $e(W)\in H_G^{d-1}(EG) = H^{d-1}(BG)$ is nonzero. Consider the natural classifying $G$-map $p : R\to EG$. In the paper~\cite{vol1992} an action of $G$ on $SO(d)$ was considered, which is equivalent to the action of $G$ on $R$, considered here. In~\cite[Proposition~3, p.~127]{vol1992} it was shown that the cohomology map $p^* : H_G^m(EG)\to H_G^m(SO(d)) = H_G^m(R)$ is injective for all $m < 2(p^k - p^{k-1})$, hence $e(s_W) = p^*e(W)\not=0$.

Thus we have $e(f)=u\times e(s_W)\not=0$ by the K\"unneth formula.

\section{Removing the limitation of strict convexity}

Consider some sequence of strictly convex smooth bodies $K_n$, that converges to $K$ in the Hausdorff metric. Each of $K_n$ has some inscribed crosspolytope $C_n$. From the compactness considerations, we can suppose that $C_n$ converges to some crosspolytope $C$, which can be degenerate (consist of one point). 

Assume that $C$ is degenerate. Again, from the compactness considerations, we may assume that the guiding cross of $C_n$ tends to some cross $X$, centered at $C\in\partial K$. The cross $X$ must have nonempty intersection with $\int K$ (from the non-angularity in Theorem~\ref{inscrcross} or smoothness in Theorem~\ref{insrccross-p}) of some positive length $\varepsilon$. It means that some line of $X$, denote it $l$, intersects $\int K$ by the open segment $\sigma$ of length at least $\varepsilon$, and $\sigma$ does not contain the center of $X$. 

Consider the guiding crosses $X_n$ of $C_n$, close enough to $X$. Their corresponding lines $l_n$ intersect $\int K_n$ by some open segments $\sigma_n$, that tend to $\sigma$. It can be easily seen, that for large enough $n$ the center of $\sigma_n$ does not coincide with the center of $X_n$, so $X_n$ cannot be a guiding cross of an inscribed crosspolytope. That is a contradiction, so $C$ is non-degenerate.

It is not clear whether the non-angularity condition can be omitted in Theorem~\ref{inscrcross}, since the crosspolytopes can degenerate in this case.

\section{The case of non-convex hypersurface}

Let us prove Theorem~\ref{insrccross-p} for arbitrary hypersurface $H$, which is an embedded sphere. 

We can smoothly deform $H$ into an ellipsoid with distinct axes, let all the deformations $H_t$ be inside some ball $B$.

Let us describe the configuration space of all crosspolytopes, similar (with positive determinant) to $C$, having center in $B$. They are parameterized by $B\times I\times SO(d)$, where $B$ is for centers, $SO(d)$ is for rotations, and $I = [a, b]$ is for homothety ($0<a<b$). 

It is clear that for small enough $a$ and large enough $b$ the crosspolytopes of sizes $\le a$ and $\ge b$ cannot be inscribed into $H_t$ for any $t\in [0, 1]$. Put $L=\partial B\times I\times SO(d) \bigcup B\times \partial I\times SO(d)$. Let $G=(Z_p)^k$ act on $SO(d)$ as above. By the K\"unneth formula $H_G^*(B\times I\times SO(d), L) = u\times v\times H_G^*(SO(d))$, where $u$ is the generator of $H^d(B, \partial B)$, $v$ is the generator of $H^1(I, \partial I)$. 

Each surface $H_t$ can be considered as the (non-degenerate everywhere) zero set of some smooth function $f_t$, we can assume that $f_t$ depends on $t$ continuously, since the homotopy $H_t$ is a restriction of some isotopy of the whole $\mathbb R^d$. Consider the map $g_t: B\times I\times SO(d) \to \mathbb R^{2d}$, which maps the triple $(x, \lambda, \rho)$ to 
$$
\left(f_t(x+\lambda\rho(e_1)), f_t(x+\lambda\rho(-e_1)), \ldots, 
f_t(x+\lambda\rho(e_d)), f_t(x+\lambda\rho(-e_d))\right).
$$
The map $g_t$ commutes with the action of $G$ on $B\times I\times SO(d)$, and the action of $G$ on $\mathbb R^{2d}$ by permuting the coordinates. Thus it makes sense to consider the Euler class $e(g_t)\in H_G^{2d}(B\times I\times SO(d), L)$. This class does not change under the deformations, now we calculate it for the case of $H_t$ being an ellipsoid. Similar to what is done in Section~\ref{obstruction} we find that in this case $e(g_t) = u\times v\times e(W)\not = 0$. Hence the zero set of $g_t$ is nonempty, every $H_t$ has an inscribed crosspolytope. Moreover, similar to the results of~\cite{vol1992}, the family of inscribed crosspolytopes of $H$ has dimension at least $\dfrac{(d-1)(d-2)}{2}$.

\section{Conclusion}

The author thanks A.Yu.~Volovikov for the discussion of these results. Some application of the above technique to the measure partition problem is given in~\cite{kar2009mp}.

\end{document}